\documentclass[nopreprintline, 1p,times,a4paper]{elsarticle}
\makeatletter
\def\ps@pprintTitle{%
 \let\@oddhead\@empty
 \let\@evenhead\@empty
 \def\@oddfoot{\centerline{\thepage}}%
 \let\@evenfoot\@oddfoot}
\makeatother

\usepackage{amssymb,amsmath}
\usepackage[mathscr]{eucal}

\journal{Linear Algebra and Its Applications}
\newtheorem{thm}{Theorem}  

\newtheorem{lemma}[thm]{Lemma}

\newproof{pf}{Proof}
\newdefinition{rmk}{Remark}

\newcommand{\sX}{\mathscr{X}}  \newcommand{\sY}{\mathscr{Y}}

\DeclareMathOperator{\Gal}{Gal}

\begin{document}

\begin{frontmatter}



\title{Nonexistence of Exceptional 5-class Association Schemes with Two $Q$-polynomial Structures   }

\author[JM]{Jianmin Ma}
\ead{Jianminma@yahoo.com}
\author[KW]{Kaishun Wang}
\ead{wangks@bnu.edu.cn}


\address[JM]{Hebei Key Lab of Computational Mathematics \& Applications \\
 and \\
College of Math \& Info. Sciences \\  Hebei Normal University, Shijiazhuang 050016, China}

\address[KW]{Sch. Math. Sci. \& Lab. Math. Com. Sys.,
Beijing Normal University, Beijing  100875, China}

\begin{abstract}

In [H. Suzuki, Association schemes with multiple $Q$-polynomial structures, J. Algebraic Combin. 7 (1998) 181-196],
Suzuki gave a classification of association schemes with multiple $Q$-polynomial structures, allowing for one exceptional
case which has five classes. In this paper, we rule out the existence of this case. Hence Suzuki's theorem   mirrors exactly the  well-known counterpart  for  association schemes with
multiple $P$-polynomial structures,  a result due to Eiichi Bannai and Etsuko Bannai in 1980.

\end{abstract}

\date{6/17/2013}
\end{frontmatter}

\section{Introduction}
\label{s:intro}

 Eiichi Bannai and Etsuko Bannai  \cite{Bannai80} studied association schemes with multiple $P$-polynomial structures and obtained  the following classification. See \cite{Bannai80} or \cite[Section III.4]{BI}.

  \begin{thm} \label{2P}
  Suppose that $\sX$ is a symmetric  association scheme
  with a $P$-polynomial structure $A_0, A_1, \dots, A_d$ for the adjacency matrices. If $\sX$ is not a polygon and has another   $P$-polynomial structure, then the new  structure is one of the following:
 \begin{enumerate}[(I)]
	\item
	$A_0, A_2, A_4, A_6, \dots, A_5, A_3, A_1;$
	\item
	$A_0, A_d, A_1,  A_{d-1},  A_2, A_{d-2}, A_3, A_{d-3}, \dots;$
	\item
	$A_0, A_d, A_2, A_{d-2}, A_4, A_{d-4}, \dots, A_{d-5}, A_5, A_{d-3}, A_3, A_{d-1}, A_1;$
	\item
	$A_0, A_{d-1}, A_2, A_{d-3}, A_4, A_{d-5}, \dots, A_5, A_{d-4}, A_3, A_{d-2}, A_1, A_d$.
	\end{enumerate}
Hence, $\sX$ admits at most two $P$-polynomial structures.
\end{thm}

 The parametric conditions for each case above can be found in \cite[Section III.4]{BI}. The question was raised whether similar result could be obtained for association schemes with multiple Q-polynomial structures in  \cite[Section III.4]{BI} and \cite{Bannai80}. Suzuki \cite{Suzuki98} settled this question  in the following theorem in 1998.

\begin{thm}  \label{suzuki}
 Suppose that $\sX $ is a symmetric  association scheme with a $Q$-polynomial structure
  $E_0, E_1, \dots, E_d$ for the primitive idempotents. If $\sX$ is not a polygon and has another   $Q$-polynomial structure, then the new structure is one of the following:
\begin{enumerate}[(I)]
	\item
	$E_0, E_2, E_4, E_6, \dots, E_5, E_3, E_1;$
	\item
	$E_0, E_d, E_1,  E_{d-1},  E_2, E_{d-2}, E_3, E_{d-3}, \dots;$
	\item
	$E_0, E_d, E_2, E_{d-2}, E_4, E_{d-4}, \dots, E_{d-5}, E_5, E_{d-3}, E_3, E_{d-1}, E_1;$
	\item
	$E_0, E_{d-1}, E_2, E_{d-3}, E_4, E_{d-5}, \dots, E_5, E_{d-4}, E_3, E_{d-2}, E_1, E_d;$
	\item
	{$d=5$ and $E_0, E_5, E_3, E_2, E_4, E_1$.}
	\end{enumerate}
Hence, $\sX$ admits at most two $Q$-polynomial structures.
\end{thm}

Note that case (V) has no counterpart  in Theorem \ref{2P}. In fact, its counterpart $A_0, A_5, A_3, A_2, A_4, A_1$ did appear in the original statement of Theorem \ref{2P} in \cite{Bannai80} but was eliminated with an easy combinatorial argument. It has been wondered if this case can also be eliminated, e.g. \cite[p.1506]{Martin09T}. We will do exactly this in this  paper, and so Theorems \ref{2P} and \ref{suzuki} are dual to each other.

\begin{thm} [Main theorem]
Case (V) in Theorem \ref{suzuki} does not occur.
\end{thm}

We conclude this section with the outline of our proof.  Any automorphism of the splitting field for an association scheme induces a permutation of its primitive idempotents. In particular, when applied to a given $Q$-polynomial structure, each nontrivial automorphism induces another  $Q$-polynomial structure. By Theorem~\ref{suzuki}, the Galois group has order at most 2. If the Galois group has order 2 and fixes the Krein parameters, we can compare the  Krein parameters from the two $Q$-polynomial structures and then determine the  Krein matrix $B_1^*$ for a putative association scheme $\sX$ in case (V). It turns out that  $B^*_1$ is completely determined by the first multiplicity $m$.   $\sX$ has a fusion scheme with $3$ classes, whose parameters can be obtained from $B_1^*$. Using elementary number theory we determine the possible values for $m$  and further show that these values give infeasible parameters of the original scheme.
If the Galois group is trivial or it has order 2 but does not fix the Krein parameters, we will derive two equations from the two $Q$-polynomial structures, which lead to a contradiction. All calculations are verified with the software package Maple 14 .

All association schemes in this paper are symmetric. The reader is referred to \cite{BI,Suzuki98b} for  missing definitions.
For recent activities on $P$- or $Q$-polynomial association schemes and related topics, see the recent survey paper by Martin and Tanaka \cite{Martin09T}.

\section{Preliminaries}
\label{prelim}
In this paper, we adopt the notation  in \cite{BI}. Let $\sX$ be a symmetric association scheme with adjacency matrices $A_i$ and primitive idempotents $E_i$, $0\le i\le d$.  Then  $A_0, \dots A_d$ span an algebra $\mathscr{M}$
 over the real field $\mathbb{R}$, called
the Bose-Mesner algebra of $\sX$.

 The intersection numbers $p^k_{ij} $
 and Krein parameters $q^k_{ij}$
are defined by
$$
 A_i A_j  = \sum_{k=0}^d p^k_{ij} A_k
 \quad \mbox{ and } \quad
E_i \circ E_j  = |X|^{-1}\sum_{k=0}^d q^k_{ij} E_k,
$$
where $\circ$ is  the entry-wise (Hadamard)  product. The intersection  matrix $B_i$ has $(j,k)$-entry $p^k_{ij}$ and the Krein matrix
$B_i^* $ has $(j,k)$-entry $q^k_{ij}$.
 When there is a possibility of confusion, we write $p^k_{i,j}$ for $p^k_{ij}$ and do the same for Krein parameters.

The adjacency matrices $A_i$ and primitive idempotents $E_i$ form two bases for $\mathscr{M}$ and so we can write down
the transition matrices between them:
$$
A_i = \sum_{j=0}^d  p_i (j)  E_j \qquad \mbox{ and } \qquad
E_i = |X|^{-1} \sum_{j=0}^d q_i(j) A_j.
$$
The numbers $p_i(j)$ and $q_i(j)$ are called eigenvalues and dual eigenvalues of $\sX$. Set
$k_i = p_i(0)$ and $m_i= q_i(0)$;  $k_i$ is  called the valency of $A_i$ and $m_i$ is the multiplicity of $E_i$.
Set $P = [p_i(j)]$,  and  $Q = [q_i(j)]$,  whose $(j,i)$-entries are $p_i(j)$ and $q_i(j)$, respectively. $P$ and $Q$ are called the first eigenmatrix and the second eigenmatrix of $\sX$.

We call $A_0, A_1, \dots, A_d$ a {\em $P$-polynomial structure} for $\sX$  if the following conditions are satisfied: for $0\le i,j,k\le d$,
\begin{enumerate}[\phantom{Item}(P1)]
\item 
$p^k_{ij} =0$ if one of  $i,j,k$ is greater than the sum of other two;
  \item 
   $p^k_{ij} \not=0$ if one of $i,j,k$ is equal to the sum of other two.
\end{enumerate}
In this case we write $c_i = p^i_{1, i-1}, a_i = p^i_{1 i} $, and  $b_i = p^i_{1, i+1}$.

Dually, we call $E_0, E_1, \dots, E_d$  a {\em $Q$-polynomial structure} for $\sX$ if, for $0\le i,j,k\le d$,
\begin{enumerate}[\phantom{Item}(Q1)]
\item [(Q1)]
$q^k_{ij} =0$ if one of  $i,j,k$ is greater than the sum of other two;
  \item[(Q2)]
   $q^k_{ij} \not=0$ if one of $i,j,k$ is equal to the sum of other two.
\end{enumerate}
Similarly, we write $c_i^* = q^i_{1, i-1}, a_i^* = q^i_{1 i} $, $b_i^* = q^i_{1, i+1}$.    Note $b_0^* = m_1$, $c^*_1 =1$, and by (Q2), $c_i^* \ne 0$ $( 1\le i\le d)$ and $b^*_i \ne 0$ $( 0\le i\le d-1)$ .

An association scheme $\sX$ is called $P$- (resp. $Q$-) polynomial if it has at least
one $P$- (resp. $Q$-) polynomial structure.
 For a $P$-polynomial scheme, each $B_i$ is a polynomial
of $B_1$ of degree $i$. Similar assertion holds for $Q$-polynomial schemes;  more precisely,
\begin{equation}
\label{B+}
B_i^*  = (B_1^* B^*_{i-1} - a_{i-1}^* B^*_{i-1} - b_{i-2}^* B^*_{i-2})/c_i^*,
\qquad  2 \le i\le d,
\end{equation}
where $B^*_0  = I$, the identity matrix of order $d+1$.

The splitting field for an association scheme  $\sX$ is the Galois extension $\boldsymbol{L}$ of the rational field $\mathbb{Q}$ by adjunction of all the eigenvalues (hence dual eigenvalues) of $\sX$ (\cite[Section II.7]{BI},  \cite{Munemasa}).  Each automorphism  of $\boldsymbol{L}/\mathbb{Q}$ applies entry-wise to all matrices  in the Bose-Mesner algebra. Let $K$ be the extension of $\mathbb{Q}$ by adjunction of all the Krein parameters. Since Krein parameters are rational functions of the eigenvalues \cite[p.65]{BI}, $K$ is an intermediate field between $\mathbb{Q}$ and $L$.  If the Galois
group $\Gal(\boldsymbol{L}/\mathbb{Q})$,  denoted also by $\Gal(\sX)$,   is non-trivial, then $\sX$ has irrational eigenvalues. For any  $\sigma\in \Gal(\sX)$, $E_0^{\sigma}, E_1^{\sigma}, \dots, E_d^{\sigma}$  are all primitive idempotents of $\mathscr{M}$ and thus a permutation of $E_0, E_1, \dots, E_d$ . It is not hard to see that $\Gal(\sX)$ acts faithfully on the primitive idempotents. In particular, if $E_0, E_1, \dots, E_d$ is a $Q$-polynomial structure, then $E_0^{\sigma}, E_1^{\sigma}, \dots, E_d^{\sigma}$ is also a $Q$-polynomial structure.

One can build new association schemes from existing ones by merging classes.
For a scheme  $\sX=(X, \{R_i\}_{i=0}^d)$,
another scheme $\sY = (X, \{S_j\}_{j=0}^e)$  is called a fusion of $\sX$ if each $S_j$ is
a union of some subset of relations $R_i$.  For any subgroup  $H\le \mbox{Gal}(\sX)$, the orbits of $H$ on the
primitive idempotents of $\sX$ give rise to a fusion scheme. The parameters of $\sY$ can be obtained from these of $\sX$ \cite{B-S}.

Now we close this section with parametric conditions for case (V) in Theorem \ref{suzuki}, which will be needed later.
\begin{lemma}
\label{param}
\emph{ \cite{Suzuki98}}
  Theorem \ref{suzuki} (V) holds if and only if $q^5_{15} =  q^5_{25} = q^5_{45} = q^5_{55} =0 \ne q^5_{35}$ and $q^5_{34} =0$.
\end{lemma}

\section{Proof of the main theorem }\label{s:krein}
\label{s:proof}
 Let  $\sX$ be an association scheme with two $Q$-polynomial structures in Theorem~\ref{suzuki} (V). In the rest of this paper,  $E_0, E_1, \dots, E_5$ is a fixed $Q$-polynomial structure for $\sX$.   Let $B_1^*$ be the first Krein matrix of $\sX$. Now by Lemma \ref{param},
we may assume
 $$
B_1^*
 =
\left[ \begin {array}{cccccc}
0&1&&&& \\\noalign{\medskip}
m& a_1^* & c_2^* &&&\\\noalign{\medskip}
& b_1^* & a_2^* & c_3^* &&\\\noalign{\medskip}
& & b_2^* &  a_3^* & c_4^* &\\\noalign{\medskip}
&& & b_3^* &  a_4^* & m\\\noalign{\medskip}
&& &&  b_4^* &  0 \end {array} \right]
:=
\left\{ \begin {array}{cccccc}
 * & 1  & c_2^*&c_3^* &c_4^* &m   \\\noalign{\medskip}
0 & a_1^*  & a_2^*  & a_3^*  &a_4^* & 0 \\\noalign{\medskip}
m  & b_1^*  & b_2^* & b_3^* &b_4^* & * \\\noalign{\medskip}
 \end {array} \right\},
$$
 where we use the notation of \cite[p.189]{BI} for tridiagonal matrices. Note the columns of $B_1^*$ sum to the rank $m (=m_1)$ of $E_1$.

We use  a hat  $\,\hat{}\,$ to indicate the second $Q$-polynomial structure.  Theorem \ref{suzuki} says
 \begin{equation}
\label{Q2}
E_{\hat{0}} = E_0, \quad E_{\hat{1}} = E_5, \quad E_{\hat{2}} = E_3, \quad E_{\hat{3}} = E_2, \quad
E_{\hat{4}}  =E_4, \quad E_{\hat{5}} = E_1.
\end{equation}
By the remarks from the previous section, $|\Gal(\sX)|\le 2$ and $\mathbb{Q} \subseteq K \subseteq L$. So we have either $\mathbb{Q} =K\ne L$ or  $K= L$.
 Now we prove the main theorem under two separate assumptions: A)    $K \ne L$; B)
   $K = L$.

  \subsection{$K\ne   L$ }
  \label{s:proof1}
In this subsection, we assume $K\ne L$ and hence $\mathbb{Q} = K$. In this case,   $\Gal(\sX)$ is generated by   an involution $\sigma$, i.e.,   $\sigma^2 $ is the identity.
 Now $E_0^{\sigma}, E_1^{\sigma}, E_2^{\sigma}, E_3^{\sigma}, E_4^{\sigma}, E_5^{\sigma}$ is the  other Q-polynomial structure and thus $E_k^\sigma = E_{\hat{k}}$. For brevity we write  $\hat{q}^r_{s t} := q^{\hat{r}}_{\hat{s}\hat{t}}.$ Since $K = \mathbb{Q}$,  $\sigma$ fixes  the Krein parameters, i.e.,
 \begin{equation}\label{krein}
 \hat{q}^r_{s t} = q^r_{st}.
 \end{equation}

Since $m_i= q^0_{ii}$, we have $m_2 = m_{\hat{2}}\ (= m_3)$ by \eqref{krein}.  It follows that
$b^*_2 = c^*_3$, since $m_i b^*_i = m_{i+1} c^*_{i+1}$. Since $ {q}^1_{11}  = \hat{q}^5_{55}$, $a^*_1=0$ by Lemma~\ref{param}. So $b_1^*=m-1$. Similarly, $a^*_3= q^3_{13} = \hat{q}^2_{52} =0$.
Since  $B_i^*$ is a polynomial in $B_1^*$ of degree $i$,
we use Eq. (\ref{B+})  to calculate $B_{2}^*$:
$$
B^*_2=  \left[ \begin {array}{cccccc}
 0&0&1&0&0&0
 \\\noalign{\medskip}
 0&{  m-1}&{ a^*_2}&{ c^*_3}&0&0
 \\\noalign{\medskip}
 {\dfrac { (m-1) m}{{ c^*_2}}}&
{\dfrac {{ (m-1)a^*_2}}{{ c^*_2}}}&
{\dfrac {{ c^*_2}(m-1)+{{ a^*_2}}^{2}+{ c^*_3}^2 -m}{{ c^*_2}}}&
{\dfrac { a^*_2{ c^*_3} }{{ c^*_2}}}&
{\dfrac {{ c^*_3}\ { c^*_4}}{{ c^*_2}}}
&0
 \\\noalign{\medskip}
0&{\dfrac {{ (m-1)c^*_3}}{{ c^*_2}}}&
{\dfrac {  a^*_2 c^*_3 } {{ c^*_2}}}&
\dfrac {{ c^*_3}^2 +  { b^*_3}{ c^*_4} -m}{{ c^*_2}}&
{\dfrac {  a^*_4  c^*_4}{{ c^*_2}}}&
{\dfrac {{ c^*_4}\ m}{{ c^*_2}}}
 \\\noalign{\medskip}
0&0&
{\dfrac {{ b^*_3}\ { b^*_2}}{{ c^*_2}}}&
{\dfrac {{ a^*_4} b^*_3}{{ c^*_2}}}&
{\dfrac {{ c^*_4}\ { b^*_3}+{{ a^*_4}}^{2}+m{ b^*_4}
-m}{{ c^*_2}}}&
{\dfrac {{ ma^*_4}}{{ c^*_2}}}
 \\\noalign{\medskip}
0&0&0&
{\dfrac {{ b^*_4}\ { b^*_3}}{{ c^*_2}}}&
{\dfrac {{ a^*_4} b^*_4 }{{ c^*_2}}}&
{\dfrac {m( b^*_4 - 1)}{{ c^*_2}}}
 \end {array}\right ].
$$
By Lemma~\ref{param}, $q^5_{25} =0$ and so $b^*_4 = 1$.

Relation (\ref{krein}) implies $B_2^* =B_{\hat{3}}^*$. Now compute $B_{\hat{3}}^*$ with \eqref{B+} and compare
$B_{\hat{3}}^*$ with $B^*_2$ to obtain the following:
\[
b_3^* = c_2^*  c_3^*, \quad  a_4^*  = 2a^*_2 = c_2^*  c_4^*, \quad
m c^*_4 = c^*_3 (m-1).
\]
Using the fact that the columns of $B^*_1$ sum to $m$, we can  obtain
$$B_1^* =
  \left\{ \begin {array}{cccccc}
 * & 1  & \dfrac{m-1}{2} & \dfrac{2m}{m+1} & \dfrac{2(m-1)}{m+1} &m   \\\noalign{\medskip}
0 & 0   & \dfrac{(m-1)^2}{2(m+1)}  & 0  & \dfrac{(m-1)^2}{m+1} & 0 \\\noalign{\medskip}
m  &  m-1   & \dfrac{2m}{m+1} & \dfrac{m(m-1)}{m+1} & 1 & * \\\noalign{\medskip}
 \end {array} \right\}.
$$

 Since   $|\mbox{Gal}(\sX)|= 2$,
 $
 \{E_0\}$,   $\{E_1,  E_5\}$,   $\{E_2, E_3\}$, $\{E_4\}$
 are the orbits of $\mbox{Gal}(\sX)$  on the primitive idempotents. By the remark before Lemma \ref{param}, these orbits give rise to a fusion scheme  $\mathcal{Y}$.
 Let
\[
T_0 =\{0\}, \quad T_1 =\{1,5\}, \quad T_2 = \{2,3\}, \quad T_3 =\{4\}.
\]
Then the Krein parameters $s^k_{ij}$  of $\mathcal{Y}$ can be calculated from  $q^r_{st}$:
\begin{equation} \label{eq:y}
s^\gamma_{ij} = \sum_{\alpha\in T_i, \ \beta \in T_j} q^\gamma_{\alpha \beta}
\end{equation}
is independent of the choice  $\gamma $ in $T_k$.

Using
(\ref{eq:y}), we can obtain the Krein  matrix $C_1^*$ of $\mathcal{Y}$:
\[
C_1^* = \left[ \begin {array}{cccc}
 0&1&0&0\\\noalign{\medskip}
 2m&0& \dfrac{m-1}{2} &2\\\noalign{\medskip}
 0&m-1&{\dfrac {{m}^{2}+6m+1}{2(m+1)}}&{
\dfrac {m-1}{4(m+1)}}\\\noalign{\medskip}0&m&{\dfrac { ( m-1 ) m
}{m+1}}& {\dfrac { ( m-1 ) ^{2}}{2(m+1)}}\end {array}
 \right].
\]
Now we can obtain the second eigenmatrix $S$ of $\mathcal{Y}$ from $C_1^*$:
\[
S=
 \left[ \begin {array}{cccc}
 1&2m&4m&{m}^{2}\\\noalign{\medskip}
 1&2&-4&1\\\noalign{\medskip}
 1&{\dfrac {{m}^{2}-10m+1+\Delta}{4(m+1)}}&
 {\dfrac { ( 5{m}^{2}-2m+5+\Delta )  ( m-1
 ) }{4 ( m+1 ) ^{2}}}&
 -{\dfrac { ( 3{m}^{2}-6m+3+\Delta ) m}{2 ( m+1 ) ^{2}}}
\\\noalign{\medskip}
1&{\dfrac {{m}^{2}+1-10m-\Delta}{4(m+1)}}&
{\dfrac { ( 5{m}^{2}-2m+5-\Delta )  ( m-1 ) }{4 ( m+1 ) ^{2}
}}&
-{\dfrac { ( 3{m}^{2}-6m+3-\Delta ) m}{2 ( m+1
 ) ^{2}}}\end {array} \right],
 \]
 where
 $
\Delta =  \sqrt {( {m}^{2}-2m+9 )  ( 9{m}^{2}-2m+1 ) }.
 $
 The first eigenmatrix of $\mathcal{Y}$ is $(m^2+ 6m+1) S^{-1}$,  and its top-most row are the  valencies of $\mathcal{Y}$:  $1, m^2$ and
\begin{equation} \label{eq:2}
\dfrac {m( 7{m}^{2}-22m+ 7)} {\sqrt { ( {m}^{2}-2m+9
 )  ( 9{m}^{2}-2m+1 ) } }  \pm 3m.
 \end{equation}

In the rest of this subsection, we determine the possible values of $m$. Since a valency is a positive integer,  the ratio
\begin{equation}\label{eq:quo}
\dfrac {m( 7{m}^{2}-22m+ 7)} {\sqrt { ( {m}^{2}-2m+9
 )  ( 9{m}^{2}-2m+1 ) } }
 \end{equation}
 is an integer, and in addition the following must hold:  (i)
 $\sqrt { ( {m}^{2}-2m+9 )  ( 9{m}^{2}-2m+1 ) } $  divides $m( 7{m}^{2}-22m+ 7)$; (ii)
$ ( {m}^{2}-2m+9 )  ( 9{m}^{2}-2m+1 ) $ is a perfect square.

 If the ratio (\ref{eq:quo}) is an integer, then any prime $p$ which divides $m^2 -2m+9$ divides either $m$ or $7m^2-22m+7$. If $p$ divides $m$ then clearly $p =3$. If $p$ divides $7m^2-22m+7$ then $p$ divides
 $$7m^2-14m+63 - (7m^2-22m+7) = 8m+56. $$
Hence either $p=2$ or $p$ divides $m+7.$ In the second case,
$p$ divides
$$(m^2-2m+9) -(m+7) = (m-1)(m-2).$$ Now $p$ divides $8$ if $m \equiv 1$ (mod p) and
$p$ divides $9$ if $m \equiv 2$ (mod p). Hence $(m-1)^2 +8$ has the form
$2^a 3^b$. The greatest common divisor of $7m^2 -22m+7$
and $9m^2-2m+1$ is also very restricted. If $p$ is an odd prime which divides
both $7m^2-22m +7$ and $9m^2-2m+1$, then (subtracting the first from the second and dividing by $2$), $p$ divides $m^2 +10m-3$. Then (subtracting the first equation from $7$ times the new equation), $p$ divides $92m-28  = 4(23m-7)$ and hence $p$ divides $7m^2 +m = 7m^2-22m +7 + (23m-7).$ Since $p$ does not divide $m$, $p$ divides $7m+1$, and thus divides $(23m-7)-3(7m+1) = 2(m-5).$ But now  $p$ also divides $7(m-5) + 36$, so $p$ divides $36$ and $p=3$.

Hence we see that $m^2 -2m+9$  has the form $2^a 3^b$ and $9m^2 - 2m+1$ has the
form $2^e 3^f$ for integers $a,b,e,f$. But now notice that the greatest common
divisor $d$ of $m^2 -2m+9$ and  $9m^2 - 2m+1$ divides $16m-80$, which is the latter equation minus 9 times the former. Also, $16$ does not divide $d$, since $m^2 -2m+9= (m-1)^2 +8$ is not divisible by 16
(and neither is $9m^2 - 2m+1= 8m^2 +(m-1)^2$). If $h$ is the odd part of $d$, then $m \equiv 5$ (mod h), and $h$ divides $24$. Hence $d$ divides $24$, since $m^2-2m+9 = (m-5)^2 + 8(m-5) + 24$. If $m$ is divisible by $3$,  then $9m^2 - 2m+1$ must be a power of $2$, but is not divisible by $16$, yielding $m=0$ under current assumptions. Hence we may suppose that $m$ is not divisible by $3$. If $9$ divides $m^2-2m +9$, then $ m \equiv 2 $ (mod 9) so $9m^2 - 2m+1 \in \{3,12,24\}$, a contradiction. Hence $m^2-2m +9  = 24$  and $m= 5$
or  $m^2-2m +9  = 8$  and $m= 1$ are the only remaining possibilities (for positive $m$ not divisible by 3, as we are currently assuming). Thus $m=1 $ and $m=5$ are the only positive integers which yield an integral value.

 Since  $c^*_2 = (m-1)/2 >0$, $m =5$. Hence $B_1^*$ is determined. Now it is a  routine exercise
to obtain the second eigenmatrix $Q$  from $B_1^*$:
\[
Q =
\left[ \begin {array}{cccccc}
1&5&10&10&25&5\\\noalign{\medskip}
1&1&-2&-2&1&1\\\noalign{\medskip}
1&
-2+\dfrac{\sqrt {21}}{3} &
\dfrac{2(1-\sqrt {21})}{3} &
 \dfrac{2(1+ \sqrt {21})}{3} &
\dfrac{5}{3} & -2- \dfrac{\sqrt {21}}{3} \\\noalign{\medskip}
1& -2- \dfrac{\sqrt {21}}{3} & \dfrac{2(1+ \sqrt {21})}{3}&
  \dfrac{2(1- \sqrt {21})}{3} &\dfrac{5}{3} &
  -2+  \dfrac{\sqrt {21}}{3} \\\noalign{\medskip}
  1&
  1+\dfrac{2\sqrt {21}}{3} & \dfrac{2(4+2\sqrt {21})}{3} &
 \dfrac{2(4- 2\sqrt {21})}{3}&-{\dfrac {25}{3}}&
   1-\dfrac{2\sqrt {21}}{3}
\\\noalign{\medskip}
1&
  1-\dfrac{2\sqrt {21}}{3} & \dfrac{2(4- 2\sqrt {21})}{3} &
 \dfrac{2(4 + 2\sqrt {21})}{3}&-{\dfrac {25}{3}}&
   1+\dfrac{2\sqrt {21}}{3}
\end {array} \right].
\]

We can obtain the intersection numbers from the matrix $Q$ by \cite[Th.3.6 (ii), p.65]{BI}:
$$
B_1 =
\left[ \begin {array}{cccccc} 0&1&0&0&0&0\\\noalign{\medskip}25&{
\dfrac {72}{7}}&10&10&{\dfrac {100}{7}}&{\dfrac {100}{7}}
\\\noalign{\medskip}0&4&5&{\dfrac {20}{3}}&{\dfrac {20}{3}}&0
\\\noalign{\medskip}0&4&{\dfrac {20}{3}}&5&0&{\dfrac {20}{3}}
\\\noalign{\medskip}0&{\dfrac {20}{7}}&\dfrac{10}{3}&0&{\dfrac {20}{7}}&{\dfrac {
25}{21}}\\\noalign{\medskip}0&{\dfrac {20}{7}}&0&\dfrac{10}{3}&{\dfrac {25}{21}}&
{\dfrac {20}{7}}\end {array} \right].
$$
Since intersection numbers must be integers, we reach a contradiction. We complete the proof of our main theorem under the assumption $K \ne  L$.

  \subsection{$K = L$} \label{s:proof2}

In this subsection, we  assume $K = L$. So  $\Gal(\sX)$ is either trivial or  $\Gal(\sX)$ has order 2 but its nonidentity element does not fixed all Krein parameters.
So relation \eqref{krein} no longer holds.

  The $(5,5)$ entry in $B^*_{\hat{5}}$  is $q_{\hat{5}\hat{5}}^{\hat{5}} =a^*_1$ and  by
Lemma \ref{param},  $a^*_1 =0$ and hence $b^*_1 = m-1$.
So we begin with
\[
B_1^* =
\left\{ \begin {array}{cccccc}
 * & 1  & c_2^*&c_3^* &c_4^* &m  \\\noalign{\medskip}
0 & 0  & a_2^*  & a_3^*  &a_4^* & 0 \\\noalign{\medskip}
m  & m -1   & b_2^* & b_3^* &b_4^* & * \\\noalign{\medskip}
 \end {array} \right\}.
\]
Now compute $B_2^*$ with (\ref{B+}). We obtain $q^5_{25} = \dfrac{m(b^*_4 -1)}{c^*_2}$ and so $b^*_4 =1$
 by Lemma  \ref{param}. Now recompute  $B_2^*, B^*_3, B^*_4$, and $ B^*_5$. With the second
 $Q$-polynomial structure, the first Krein matrix ${B}^*_{\hat{1}} = (q_{\hat{1}\hat{i}}^{\hat{k}})$ with $(i,k)$ entry
 $q_{\hat{1}\hat{i}}^{\hat{k}}$  is tridiagonal with nonzero off-diagonal entries. The (2,1) entry of ${B}^*_{\hat{1}}$  is
 $\dfrac{m a^*_4}{c_2^* c_3^*}$, so $a_4^* \ne 0$.  We also have
$$
q^{\hat{4}}_{\hat{1}\hat{1}} = \dfrac{  c^*_4  b^*_3 +  {a^*_4}^2  - a^*_2  a^*_4 - (m-1) c^*_2- a^*_3 a^*_4}{ c^*_4  c^*_3  c^*_2}
\quad \mbox{ and
  }\quad
q^{\hat{4}}_{\hat{1}\hat{2}} = \dfrac{  c^*_4  b^*_3 +  {a^*_4}^2  - a^*_2* a^*_4 - (m-1) c^*_2}{ c^*_3 c^*_2},
$$
which should be 0 by (Q2). It follows that  $a^*_3 =0$.
 Now, we have
\[
B_1^* =
\left\{ \begin {array}{cccccc}
 * &1  & c_2^*&c_3^* &c_4^* &m  \\\noalign{\medskip}
0 & 0 & a_2^*  & 0  &a_4^* & 0 \\\noalign{\medskip}
m  & m-1  & b_2^* & b_3^* & 1 & * \\\noalign{\medskip}
 \end {array} \right\}.
\]
 Similarly, we compute  $B_5^*$, from which we can obtain
${B}^*_{\hat{1}}$.   We find
$$
\hat{q}^1_{11} =
\dfrac{ a_4^* m-a_2^* c_4^* b_3^* - b_2^* a_4^*c_3^* }{c_4^* c_3^*c_2^*}.
$$
By Lemma \ref{param}, $ q^5_{55} = 0$.  As $\hat{q}^1_{11} =p^5_{55}$,  we have
\begin{align}
a_4^* m-a_2^* c_4^* b_3^* - b_2^* a_4^*c_3^* & =0. \label{eq1}
\end{align}

Now we derive another  equation.
Let $v^*_i(x)~( 2\le i\le 5) $  be the polynomials defined by $B_1^*$ \cite[Section III.1]{BI}:
$$v^*_0(x) = 1, \quad v^*_1(x) = x, \quad
x v^*_i(x) = b_{i-1}^* v^*_{i-1}(x) + a_i^* v^*_i(x) + c_{i+1}^* v^*_{i+1}(x),
$$
 where $c_6^* :=1$.
So $v^*_6(x)$ annihilates $B_1^*$. Since $m$ is an eigenvalue of $B^*_1$,  we have
$$
v^*_6(m) = \frac{ m^2 ( m-1 )  ( -m^2{a_4^*}+ m{a_4^*} {c_2^*} -m{b_3^*} {c_4^*}+
  m{a_2^*} {a_4^*} + {a_2^*}{b_3^*} {c_4^*} + {a_4^*}{b_2^*}{c_3^*}+
 {c_4^*} {b_3^*} {c_2^*})}{c^*_2 c^*_3 c^*_4} =0,
$$
which gives
$$
-m^2{a_4^*}+ m{a_4^*} {c_2^*} -m{b_3^*} {c_4^*}+
  m{a_2^*} {a_4^*} + {a_2^*}{b_3^*} {c_4^*} + {a_4^*}{b_2^*}{c_3^*}+
 {c_4^*} {b_3^*} {c_2^*} =0.
$$
Using Eq. (\ref{eq1}), it simplifies to
$$
 -m^2{a_4^*}+ m{a_4^*} {c_2^*} -m{b_3^*} {c_4^*}+
  m{a_2^*} {a_4^*} + m a_4^* +
 {c_4^*} {b_3^*} {c_2^*} =0.
$$
By $a_2^* + c_2^* = m - b_2^*$, it becomes
\begin{align} \label{eqA}
m a_4^* (1 - b_2^*) = b_3^* c_4^* (m - c_2^*).
\end{align}
Taking the difference of  Eq. (\ref{eq1}) and (\ref{eqA}), we obtain
$$
m a_4^* b_2^* = b_3^* c_4^* ( a_2^* + c_2^* - m) + a_4^* b_2^* c_3^*.
$$
This simplifies to $a_4^*b_2^* b_3^* = - b_2^* b_3^* c_4^*$, as $c_3^* + b_3^* =m$ . Hence, $a^*_4 + c^*_4 = 0$, absurd.
  We have proved the main theorem under the assumption $K = L$.

\section{Concluding Remarks}

\textbf{Remark 1.}
$P$-polynomial structures (I) and (II) come in pair. As pointed out in \cite[p.243]{BI}, (I) and (II) are dual to each other, i.e., exchanging the roles of the first and second structures, the first is of type (II) (resp. (I)) in terms of the second if the second is of type (I) (resp. (II)) in terms of the first. The structures (III), (IV) are self dual. The same comment applies to $Q$-polynomial structures.

\textbf{Remark 2.} Cases (I), (II) and (IV) of Theorem \ref{2P} are realized
\cite[Section III.4]{BI}.  Similarly for Theorem \ref{suzuki}, case (I) is realized by the half cube $\frac{1}{2}H(2d+1,2)$, case (II) by the antipodal quotient of the cube $\widetilde{H}(2d+1,2)$ and the dual polar graph on $[{^2A}_{2d-1}(r)], r\ge 2$, and case (IV) by the cube $H(d,2)$, $d$ even. Note each member of these four families is a distance regular graph admitting two $Q$-polynomial structures.   Dickie \cite{Dickie} classified such graphs with diameter $d\ge 5$ and valency $k\ge 3$.
In view of Remark 1, $\frac{1}{2}H(2d+1,2)$ admits a $Q$-polynomial structure of type (II);
  $\widetilde{H}(2d+1,2)$ and $[{^2A}_{2d-1}(r)], r\ge 2,$ admit a $Q$-polynomial structure of type (I).

  Suzuki \cite{Suzuki96} showed that the class number $d\le 4$ for Th.\ref{2P} (III), and it is yet unknown whether similar result can be obtained for  Th.\ref{suzuki} (III).  We do not know any association scheme with $d\ (\ge 3)$ classes,  which admits two  $Q$-polynomial structures and no $P$-polynomial structures.

\textbf{Remark 3.} A careful reader will notice that the proof in Subsection \ref{s:proof2} does not even use the assumption $K=L$, so Subsection \ref{s:proof1} can be omitted. However, we feel that Subsection \ref{s:proof1} serves to illustrate a powerful use of   the Galois group. Plus, the elementary number theoretic proof in this subsection is very beautiful.

\section*{Acknowledgments}

We are indebted to the anonymous reviewers whose detailed and thoughtful reports save us a few mistakes and improve the readability of this paper.
We thank Hajime Tanaka  for bringing to our attention the problem of this paper during
the 2010 Workshop on Schemes and Spheres at Worcester Polytechnic Institute, and thus the first author thanks Bill Martin for his hospitality and many inspiring conversations. Geoff Robinson provides the elementary approach in Section \ref{s:proof} which determines the values of $m$ from (\ref{eq:quo}).
We also thank Eiichi Bannai for sending us his original paper \cite{Bannai80} that helps to correct a few inaccurate statements in an early draft.

 JM is partially supported by the Natural Science Foundation of Hebei province (A2012205079) and Science Foundation of  Hebei Normal University (L2011B02). KW  is partially supported by  NSF of China (11271047, 11371204) and the Fundamental Research Funds for the Central Universities of China.

\end{document}